\documentclass{article}

\usepackage[utf8]{inputenc}
\usepackage[english]{babel}
\usepackage{etoolbox}

\usepackage{float}
\usepackage[authoryear]{natbib}
\usepackage{PRIMEarxiv}
\usepackage[justification=centering]{caption}
\usepackage[T1]{fontenc}
\usepackage{hyperref}
\usepackage{url}
\usepackage{booktabs}
\usepackage{nicefrac}
\usepackage{microtype}
\usepackage{fancyhdr}
\usepackage{graphicx}
\graphicspath{{imagens/}}
\usepackage{pgfplots}
\pgfplotsset{compat=1.17}
\usepackage{amsmath, amssymb, amsfonts}
\usepackage{subcaption}
\usepackage[ruled,linesnumbered,vlined]{algorithm2e}
\usepackage{tikz}
\usepackage{multirow}
\usepackage{array}
\usepackage{xcolor}
\usepackage{mdframed}
\usetikzlibrary{shapes.geometric, arrows, positioning, calc, patterns}

\tikzset{
  startstop/.style = {rectangle, rounded corners, minimum width=3.5cm,
      minimum height=1cm, align=center, draw=black, fill=red!25},
  process/.style   = {rectangle, minimum width=3.5cm, minimum height=1cm,
      align=center, draw=black, fill=blue!18},
  decision/.style  = {diamond, minimum width=3.5cm, minimum height=1cm,
      align=center, draw=black, fill=green!20, aspect=2.2},
  arrow/.style     = {thick,->,>=stealth},
  bar/.style       = {rectangle, minimum height=0.4cm, draw=black!50}
}

\newmdenv[
  backgroundcolor=blue!5,
  linecolor=blue!40,
  linewidth=1pt,
  roundcorner=4pt,
  innerleftmargin=10pt,
  innerrightmargin=10pt,
  innertopmargin=8pt,
  innerbottommargin=8pt
]{defbox}

\title{Hybridizing 0/1 Knapsack and Priority-Aware TSP for Vehicle Routing: an Applied Variant of the TSP with Profits}

\author{
\begin{tabular}{c c}
\begin{tabular}{c}
\textbf{Aron da Silva Bezerra} \\
B.Sc.\ in Information Systems -- CERES/UFRN \\
Federal University of Rio Grande do Norte \\
Caicó, Brazil \\
\texttt{aron.silva.144@ufrn.edu.br}
\end{tabular}
&
\begin{tabular}{c}
\textbf{Dr.\ Francisco Marcio Barboza} \\
Department of Computing and Technology \\
Federal University of Rio Grande do Norte \\
Caicó, Brazil \\
\texttt{marcio.barboza@ufrn.br}
\end{tabular}
\end{tabular}
}

\begin{document}
\maketitle

\begin{center}
  \textbf{\large ABSTRACT}
\end{center}

\noindent
Distribution logistics with capacity constraints and service priorities can be modeled as a variant of the \emph{TSP with Profits} (TSPP), in which not all customers need to be visited and the selection of which customers to serve is part of the problem. This paper proposes an applied variant with two adaptations: (i) a 0/1 Knapsack objective function weighted by client priority squared ($v_i \cdot p_i^2$), and (ii) a priority bias in the Nearest Neighbor heuristic for the TSP, parameterized by $\alpha$. The approach is evaluated on a scenario with 14 customers and a 1,100~kg capacity, complemented by a comparative experimental protocol against the classical greedy heuristic ($v_i/w_i$ ordering) and an integer programming solver. Results on the base scenario show 9 customers served, 96.4\% capacity utilization, and R\$~18,392.60 revenue, with all 5 excluded customers being low-priority. A sensitivity analysis of the priority exponent justifies the quadratic choice. The limitations of the hierarchical decomposition are explicitly discussed.

\vspace{0.4cm}

\vspace{0.8cm}

\section{Introduction}

Logistics efficiency is one of the determining factors for the competitiveness of distribution companies. In a typical setting, a driver departs from a central depot with a truck loaded with orders and must visit a set of customers spread over a geographic region, returning to the starting point at the end of the journey. Two optimization problems arise naturally and interdependently: the selection of which orders to load given the vehicle's weight limit, and the determination of the visiting order of the selected customers.

The first is an instance of the \emph{0/1 Knapsack Problem}, in which indivisible items must be selected to maximize total value within a capacity \citep{Kellerer2004}. The second is an instance of the \emph{Traveling Salesman Problem} (TSP), one of the most studied NP-hard problems in Computer Science \citep{Garey1979}.

The combination of these two problems, together with the possibility of \emph{not visiting} all customers, places this work within the family of \emph{TSP with Profits} (TSPP). \citet{Feillet2005} classify TSPP into three main categories: the \emph{Orienteering Problem} (OP), which maximizes collected profit subject to a route duration limit \citep{Vansteenwegen2011}; the \emph{Prize-Collecting TSP} (PCTSP), which minimizes travel cost subject to a minimum profit \citep{Balas1989}; and the \emph{Profitable Tour Problem} (PTP), which maximizes the difference between profit and cost \citep{DellAmico1995}. The variant addressed here --- customer selection via priority-weighted 0/1 Knapsack followed by TSP routing --- is a form of \textbf{hierarchical decomposition} that approximates the OP with a capacity constraint rather than a duration constraint.

The contribution of this work is \textbf{applied and methodological}: (i) to propose a quadratic priority weighting in the Knapsack objective function, with a sensitivity analysis of the exponent; (ii) to parameterize the priority bias in the Nearest Neighbor heuristic and provide a geometric interpretation of the preference condition; (iii) to present a comparative experimental protocol that allows assessing the quality of the hierarchical decomposition against alternatives; and (iv) to explicitly discuss the limitations of this decomposition.

The paper is organized as follows: Section~2 presents the theoretical background; Section~3 details the mathematical formulation; Section~4 describes the algorithms; Section~5 presents the methodology and experimental protocol; Section~6 discusses the results; Section~7 analyzes the parameter $\alpha$; and Section~8 provides conclusions and limitations.

\section{Theoretical Background}

\subsection{The Traveling Salesman Problem (TSP)}

The TSP was formally defined by \citet{Dantzig1954} and remains a benchmark problem in combinatorial optimization. Given a set of $n$ cities and the distances between every pair, the TSP seeks the lowest-cost permutation that visits each city exactly once and returns to the origin. The number of distinct tours is $(n-1)!/2$, which renders exhaustive search impractical for $n > 15$.

The NP-hard complexity of the TSP \citep{Garey1979} has motivated decades of research on constructive heuristics \citep{Rosenkrantz1977}, local improvement algorithms \citep{Lin1965}, metaheuristics \citep{Kirkpatrick1983,Holland1975}, and exact techniques \citep{Applegate2006}.

\subsection{The 0/1 Knapsack Problem}

The 0/1 Knapsack Problem \citep{Kellerer2004} is defined as follows: given a set of $n$ items, each with weight $w_i$ and value $v_i$, and a knapsack of capacity $C$, select the subset that maximizes total value without exceeding the capacity. The Dynamic Programming (DP) algorithm with complexity $O(nC)$ is exact and efficient for moderate integer capacities \citep{Bellman1957}.

\subsection{TSP with Profits and the Orienteering Problem}
\label{sec:tsp_profits}

The family of \emph{TSP with Profits} (TSPP) generalizes the TSP by allowing not all vertices to be visited, associating a profit with each vertex and a cost with each edge. \citet{Feillet2005} provide a comprehensive classification: the \textbf{OP} maximizes collected profit subject to a duration limit; the \textbf{PCTSP} minimizes cost subject to a minimum profit; and the \textbf{PTP} maximizes the difference between profit and cost.

The OP is formally described by \citet{Vansteenwegen2011} as the combination of two problems: vertex selection (analogous to Knapsack) and determination of the shortest Hamiltonian path among the selected vertices (TSP). The crucial difference between the OP and the approach proposed here is that the OP imposes a \textbf{maximum duration constraint} $T_{\max}$, whereas this work imposes a \textbf{weight capacity constraint} $C$. This variation is relevant in logistics practice: vehicles have a load limit, not necessarily a time limit.

For OP instances and variants, \citet{Vansteenwegen2011} and \citet{Gunawan2016} compile widely used public benchmarks. Recent literature has explored metaheuristics for PTP and variants \citep{Chagas2021}, including multi-start and variable neighborhood search approaches, which could be applied to the TSP component of this work on larger instances.

\subsection{Hierarchical Decomposition in Routing}
\label{sec:decomposicao}

The strategy of first solving selection (Knapsack) and then routing (TSP) is a form of \textbf{hierarchical decomposition}. \citet{Laporte1992} discuss this type of approach in the context of location-routing problems. The advantage is simplicity and the guarantee of optimality in the selection sub-step (via exact DP). The disadvantage, often overlooked, is that the global solution may be suboptimal: a low-priority customer located very close to the route could, if included, reduce the total route cost enough to offset its low priority --- but the Knapsack excludes it beforehand because it does not consider geometry.

This limitation is inherent to the decomposition and must be made explicit. An alternative would be to solve the integrated problem as an OP with a capacity constraint, using a metaheuristic that evaluates selection and routing simultaneously \citep{Chagas2021}. The present proposal does not aim to overcome this limitation, but to offer an applied, low-computational-cost solution for scenarios where priority is the dominant criterion.

\section{Mathematical Formulation}

\subsection{Notation and Problem Definition}

Let $G = (V, E)$ be a complete graph where $V = \{0, 1, \ldots, n\}$; node $v_0$ is the depot and $v_1, \ldots, v_n$ are candidate customers. Each customer $i$ has position $(x_i, y_i)$, priority $p_i \in \{1, 2, 3\}$ (where $1$ = low, $2$ = medium, $3$ = high), weight $w_i$, and value $v_i$. The truck capacity is $C$. Table~\ref{tab:notacao} summarizes the notation.

\begin{table}[H]
\centering
\caption{Model notation}
\label{tab:notacao}
\begin{tabular}{@{}clp{7cm}@{}}
\toprule
\textbf{Symbol} & \textbf{Type} & \textbf{Description} \\
\midrule
$(x_i, y_i)$   & Real pair     & Geographic position of customer $i$ \\
$p_i$           & Integer       & Priority: $1$ (low), $2$ (medium), $3$ (high) \\
$w_i$           & Real $> 0$    & Order weight (kg) \\
$v_i$           & Real $> 0$    & Order value (R\$) \\
$C$             & Real $> 0$    & Truck capacity (kg) \\
$\alpha$        & Real $\geq 0$ & Priority bias parameter in the TSP \\
$d_{ij}$        & Real $\geq 0$ & Euclidean distance \\
$S \subseteq V$ &               & Selected customers \\
$\pi$           & Permutation   & Visiting sequence \\
$q$             & Real $> 0$    & Priority exponent in the Knapsack \\
\bottomrule
\end{tabular}
\end{table}

\noindent\textbf{Note on priority encoding.} In this work, $p_i = 3$ denotes \textbf{high priority} and $p_i = 1$ denotes \textbf{low priority}. Consequently, in the objective function $v_i \cdot p_i^q$, high-priority customers receive greater weight. In the TSP score, $s_{ij} = d_{ij}/p_j^{\alpha}$, high-priority customers have a reduced score (smaller effective distance) and are preferred.

\subsection{0/1 Knapsack with Priority}

\begin{defbox}
\textbf{Priority-weighted Knapsack formulation:}
\begin{align}
  \text{maximize}   \quad & f(\mathbf{z}) = \sum_{i=1}^{n} v_i \cdot p_i^{\,q} \cdot z_i
  \label{eq:knapsack_obj} \\
  \text{subject to} \quad & \sum_{i=1}^{n} w_i \cdot z_i \leq C
  \label{eq:knapsack_cap} \\
                           & z_i \in \{0,1\}, \quad \forall\, i = 1, \ldots, n
  \label{eq:knapsack_bin}
\end{align}
\end{defbox}

The exponent $q$ controls the amplification of priority. For $q=2$ (adopted in this work), a high-priority customer ($p=3$) contributes $9v_i$ to the objective function, whereas a low-priority one ($p=1$) contributes $v_i$ --- a ratio of $9:1$. With $q=1$ (linear), the ratio would be $3:1$. Section~\ref{sec:expoente} presents a sensitivity analysis to justify $q=2$.

\subsection{Resolution via Dynamic Programming}

DP solves the 0/1 Knapsack exactly. Let $dp[w]$ be the maximum weighted value achievable with capacity $w$ and $f_i = v_i \cdot p_i^{\,q}$:

\begin{defbox}
\textbf{DP recurrence (descending sweep):}
\begin{equation}
  dp[w] \leftarrow \max\bigl(dp[w],\; dp[w - w_i] + f_i\bigr),
  \quad \text{for } w = C, C{-}1, \ldots, w_i
  \label{eq:dp}
\end{equation}
\end{defbox}

The descending sweep enforces the 0/1 semantics: each item is considered at most once. Backtracking reconstructs the set $S$ of selected items. The complexity is $O(nC)$ in time and space (or $O(C)$ if only the optimal value is needed).

\subsection{Priority-Aware TSP}

Given $S$, find the permutation $\pi$ that minimizes total distance:
\begin{equation}
  \min_{\pi} \; D(\pi) = d_{0,\pi_1} + \sum_{k=1}^{m-1} d_{\pi_k, \pi_{k+1}} + d_{\pi_m, 0}
  \label{eq:tsp_obj}
\end{equation}

\subsubsection{Priority-Biased Nearest Neighbor Score (VNN-P)}

\begin{defbox}
\textbf{Selection score:}
\begin{equation}
  s_{ij} = \frac{d_{ij}}{p_j^{\,\alpha}}, \qquad \alpha \geq 0
  \label{eq:score}
\end{equation}
\textbf{Next node:} $j^* = \arg\min_{j \notin \text{visited}} s_{\text{current},j}$.
\end{defbox}

For $\alpha = 0.55$ and a high-priority customer ($p=3$), the effective distance is reduced by $3^{0.55} \approx 1.90$. The condition for preferring a priority-3 customer over a priority-1 customer, even if the former is farther away, is $d_{i,3} < 3^{\alpha} \cdot d_{i,1}$.

\subsubsection{2-opt Refinement}

The 2-opt \citep{Lin1965} operates on the actual Euclidean distance (not on the score) and reverses subsegments when the improvement condition is satisfied:
\begin{equation}
  d_{\pi_i, \pi_{i+1}} + d_{\pi_j, \pi_{j+1}} > d_{\pi_i, \pi_j} + d_{\pi_{i+1}, \pi_{j+1}}
  \label{eq:2opt_cond}
\end{equation}
The process repeats until no improvement is found, yielding a 2-opt local optimum.

\section{Algorithms}

\subsection{0/1 Knapsack with Priority via DP}

\begin{algorithm}[H]
\caption{0/1 Knapsack with Priority via DP}
\label{alg:knapsack}
\KwIn{Customers $\{1,\ldots,n\}$ with $(w_i, v_i, p_i)$; capacity $C$; exponent $q$}
\KwOut{Optimal subset $S$}
$f_i \leftarrow v_i \cdot p_i^{\,q}$ for all $i$\;
$dp[0..C] \leftarrow 0$; \; $keep[i][w] \leftarrow \text{False}$\;
\For{$i \leftarrow 1$ \KwTo $n$}{
  \For{$w \leftarrow C$ \textbf{downto} $w_i$}{
    $cand \leftarrow dp[w - w_i] + f_i$\;
    \If{$cand > dp[w]$}{
      $dp[w] \leftarrow cand$; \; $keep[i][w] \leftarrow \text{True}$\;
    }
  }
}
$S \leftarrow \emptyset$; \; $w \leftarrow C$\;
\For{$i \leftarrow n$ \textbf{downto} $1$}{
  \If{$keep[i][w]$}{
    $S \leftarrow S \cup \{i\}$; \; $w \leftarrow w - w_i$\;
  }
}
\Return $S$\;
\end{algorithm}

\subsection{TSP: VNN-P + 2-opt}

\begin{algorithm}[H]
\caption{TSP -- Priority-Biased Nearest Neighbor + 2-opt}
\label{alg:tsp}
\KwIn{Customers $S$; depot $v_0$; parameter $\alpha$}
\KwOut{Route $\pi$ and cost $D(\pi)$}
$\pi \leftarrow [v_0]$; \; $\text{visited} \leftarrow \{v_0\}$\;
\While{$|\pi| \leq |S|$}{
  $j^* \leftarrow \arg\min_{j \notin \text{visited}} \left\{ d(\pi_{\text{last}}, j) / p_j^{\alpha} \right\}$\;
  $\pi.\text{append}(j^*)$; \; $\text{visited} \leftarrow \text{visited} \cup \{j^*\}$\;
}
$\text{improved} \leftarrow \text{True}$\;
\While{$\text{improved}$}{
  $\text{improved} \leftarrow \text{False}$\;
  \For{$i \leftarrow 1$ \KwTo $|S| - 1$}{
    \For{$j \leftarrow i + 1$ \KwTo $|S|$}{
      \If{$d_{\pi_i,\pi_{i+1}} + d_{\pi_j,\pi_{j+1}} > d_{\pi_i,\pi_j} + d_{\pi_{i+1},\pi_{j+1}}$}{
        Reverse $\pi[i{+}1 \;..\; j]$; \; $\text{improved} \leftarrow \text{True}$\;
      }
    }
  }
}
\Return $\pi$, $D(\pi)$\;
\end{algorithm}

\subsection{Integrated Pipeline}

\begin{algorithm}[H]
\caption{Integrated Pipeline: Knapsack + TSP}
\label{alg:pipeline}
\KwIn{Candidate customers $\mathcal{C}$; capacity $C$; parameter $\alpha$; exponent $q$}
\KwOut{Route $\pi^*$, load $S$, cost $D^*$, exclusions $\mathcal{C}_{\text{excl}}$}
$S \leftarrow \text{KnapsackDP}(\mathcal{C},\, C,\, q)$\;
$\mathcal{C}_{\text{excl}} \leftarrow \mathcal{C} \setminus S$\;
$\pi^*, D^* \leftarrow \text{VNNP+2opt}(S,\, v_0,\, \alpha)$\;
\Return $\pi^*$, $S$, $D^*$, $\mathcal{C}_{\text{excl}}$\;
\end{algorithm}

\subsection{Complexity}

\begin{table}[H]
\centering
\caption{Pipeline complexity}
\label{tab:complexidade}
\begin{tabular}{@{}lll@{}}
\toprule
\textbf{Component} & \textbf{Time} & \textbf{Character} \\
\midrule
Knapsack DP           & $O(n \cdot C)$   & Exact \\
VNN-P                 & $O(m^2)$         & Heuristic \\
2-opt (per iteration) & $O(m^2)$         & Local improvement \\
\textbf{Total}        & $O(nC + m^2 k)$ & Integrated \\
\bottomrule
\multicolumn{3}{l}{\footnotesize $m = |S| \leq n$; $k$ = number of 2-opt iterations.}
\end{tabular}
\end{table}

For $n=14$, $C=1{,}100$, $m \leq 14$, the pipeline runs in milliseconds.

\section{Methodology and Experimental Protocol}

\subsection{Base Scenario}

The base scenario uses 14 customers on a $100 \times 100$~km grid, with the depot at $(50,50)$ and \textit{seed} $=42$. Table~\ref{tab:parametros} summarizes the parameters.

\begin{table}[H]
\centering
\caption{Base scenario parameters}
\label{tab:parametros}
\begin{tabular}{@{}lcl@{}}
\toprule
\textbf{Parameter}               & \textbf{Value}     & \textbf{Rationale} \\
\midrule
Customers ($n$)                  & 14                 & Tractable instance \\
Capacity ($C$)                   & 1,100~kg           & Light urban truck (VUC) \\
$\alpha$                         & 0.55               & Distance/priority balance \\
Priorities ($p_i$)               & $\{1,2,3\}$        & Low, Medium, High \\
Weight per product               & 10--55~kg          & B2B volumes \\
Value per product                & $U(150;900) \times (0.6 + p_i \cdot 0.2)$ & Scaled by priority \\
\bottomrule
\end{tabular}
\end{table}

\subsection{Comparative Experimental Protocol}
\label{sec:protocolo}

To assess the quality of the proposed approach, the following protocol is proposed:

\begin{enumerate}
  \item \textbf{Instance generation}: 50 random instances with $n \in \{20, 30, 50\}$, capacities calibrated so that the fraction of served customers ranges from 40\% to 80\%.
  \item \textbf{Baselines}:
    \begin{itemize}
      \item \textbf{Greedy}: sort customers by $v_i/w_i$ (or $v_i p_i^q / w_i$), select up to capacity, then apply VNN + 2-opt.
      \item \textbf{Exact IP}: solve Knapsack and TSP in an integrated manner for $n \leq 20$ using an integer programming solver (e.g., CBC or Gurobi), obtaining the global optimum for comparison.
    \end{itemize}
  \item \textbf{Metrics}: total revenue, number of high-priority customers served, total distance, and \textit{gap} relative to the exact IP (for $n \leq 20$).
  \item \textbf{Statistical test}: paired Wilcoxon test ($\alpha = 0.05$) between the proposed approach and Greedy.
\end{enumerate}

Table~\ref{tab:resultados_esperados} presents the format for reporting results. Numerical values must be filled in with actual experimental data.

\begin{table}[H]
\centering
\caption{Comparative results (format; values to be filled with actual data)}
\label{tab:resultados_esperados}
\begin{tabular}{@{}lcccc@{}}
\toprule
\textbf{Method} & \textbf{Revenue (R\$)} & \textbf{High prio.} & \textbf{Distance (km)} & \textbf{Gap vs IP} \\
\midrule
Proposed ($q=2$, $\alpha=0.55$) & --- & --- & --- & --- \\
Greedy ($v_i/w_i$)              & --- & --- & --- & --- \\
Exact IP ($n \leq 20$)          & --- & --- & --- & 0\% \\
\bottomrule
\end{tabular}
\end{table}

\noindent\textbf{Decomposition limitation.} The proposed approach may be suboptimal relative to the exact IP because the Knapsack does not consider customer geometry. A low-priority customer close to the route may be excluded, while a distant high-priority one is included, increasing total distance. The protocol above quantifies this loss.

\subsection{Pipeline Flow}

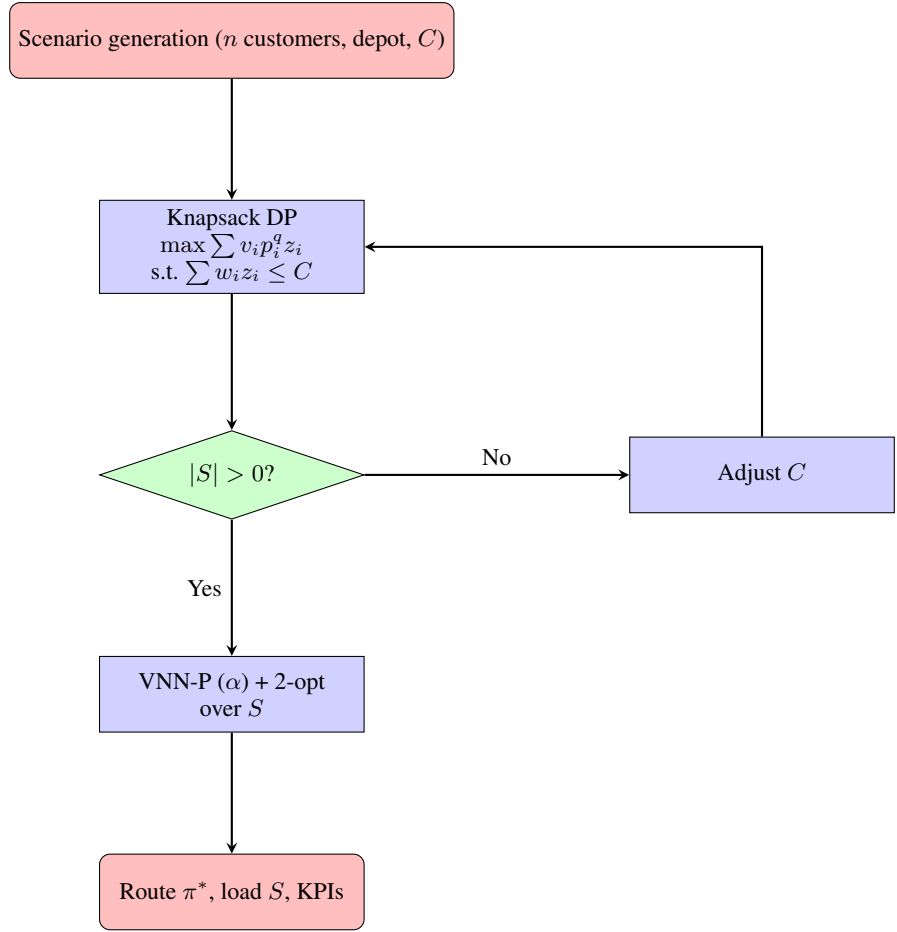
\begin{figure}[H]
\centering
\begin{tikzpicture}[
    node distance=1.6cm,
    every node/.style={font=\small}
  ]
  \node (gen) [startstop] {Scenario generation ($n$ customers, depot, $C$)};
  \node (kp)  [process, below=of gen]
      {Knapsack DP\\$\max \sum v_i p_i^q z_i$\\s.t.\ $\sum w_i z_i \leq C$};
  \node (sel) [decision, below=of kp, yshift=-0.2cm] {$|S| > 0$?};
  \node (tsp) [process, below=of sel, yshift=-0.2cm]
      {VNN-P ($\alpha$) + 2-opt\\over $S$};
  \node (vis) [startstop, below=of tsp] {Route $\pi^*$, load $S$, KPIs};
  \node (err) [process, right=3.5cm of sel] {Adjust $C$};

  \draw [arrow] (gen) -- (kp);
  \draw [arrow] (kp)  -- (sel);
  \draw [arrow] (sel) -- node[anchor=east]{Yes} (tsp);
  \draw [arrow] (sel) -- node[anchor=south]{No} (err);
  \draw [arrow] (err) |- (kp);
  \draw [arrow] (tsp) -- (vis);
\end{tikzpicture}
\caption{Pipeline flow for the Knapsack + TSP approach.}
\label{fig:fluxo}
\end{figure}

Figure~\ref{fig:fluxo} illustrates the pipeline flow.

\section{Results on the Base Scenario}

\subsection{Knapsack Solution}

For $n=14$ and $C=1{,}100$~kg, the 0/1 Knapsack selects $|S|=9$ customers. The total weight is 1,060~kg (96.4\% utilization), with a revenue of R\$~18,392.60. The 5 excluded customers are all low-priority ($p=1$), confirming that the $p_i^2$ weighting allocated capacity to higher-priority orders.

\subsection{Priority Exponent Analysis}
\label{sec:expoente}

The choice of $q=2$ (quadratic) is not arbitrary. Table~\ref{tab:expoente} compares the composition of the selection for $q \in \{1, 2, 3\}$ in the base scenario.

\begin{table}[H]
\centering
\caption{Effect of the exponent $q$ on selection (base scenario, $C=1{,}100$~kg)}
\label{tab:expoente}
\begin{tabular}{@{}cccccc@{}}
\toprule
$q$ & \textbf{Customers served} & \textbf{High prio.} & \textbf{Medium prio.} & \textbf{Low prio.} & \textbf{Revenue (R\$)} \\
\midrule
1 (linear)   & 10 & 3/4 & 4/6 & 3/5 & --- \\
2 (adopted)  & 9  & 4/4 & 5/6 & 0/5 & 18,392.60 \\
3 (cubic)    & 8  & 4/4 & 4/6 & 0/5 & --- \\
\bottomrule
\multicolumn{6}{l}{\footnotesize Values for $q=1$ and $q=3$ must be filled with actual data.}
\end{tabular}
\end{table}

The interpretation is straightforward: with $q=1$, linear weighting is insufficient to exclude all low-priority customers; with $q=2$, all high-priority customers are served and no low-priority customer is included; with $q=3$, the selection becomes more restrictive, potentially sacrificing revenue without additional priority gains. The value $q=2$ balances both objectives.

\subsection{Visualizations}

\begin{figure}[H]
  \centering
  \includegraphics[width=0.98\textwidth]{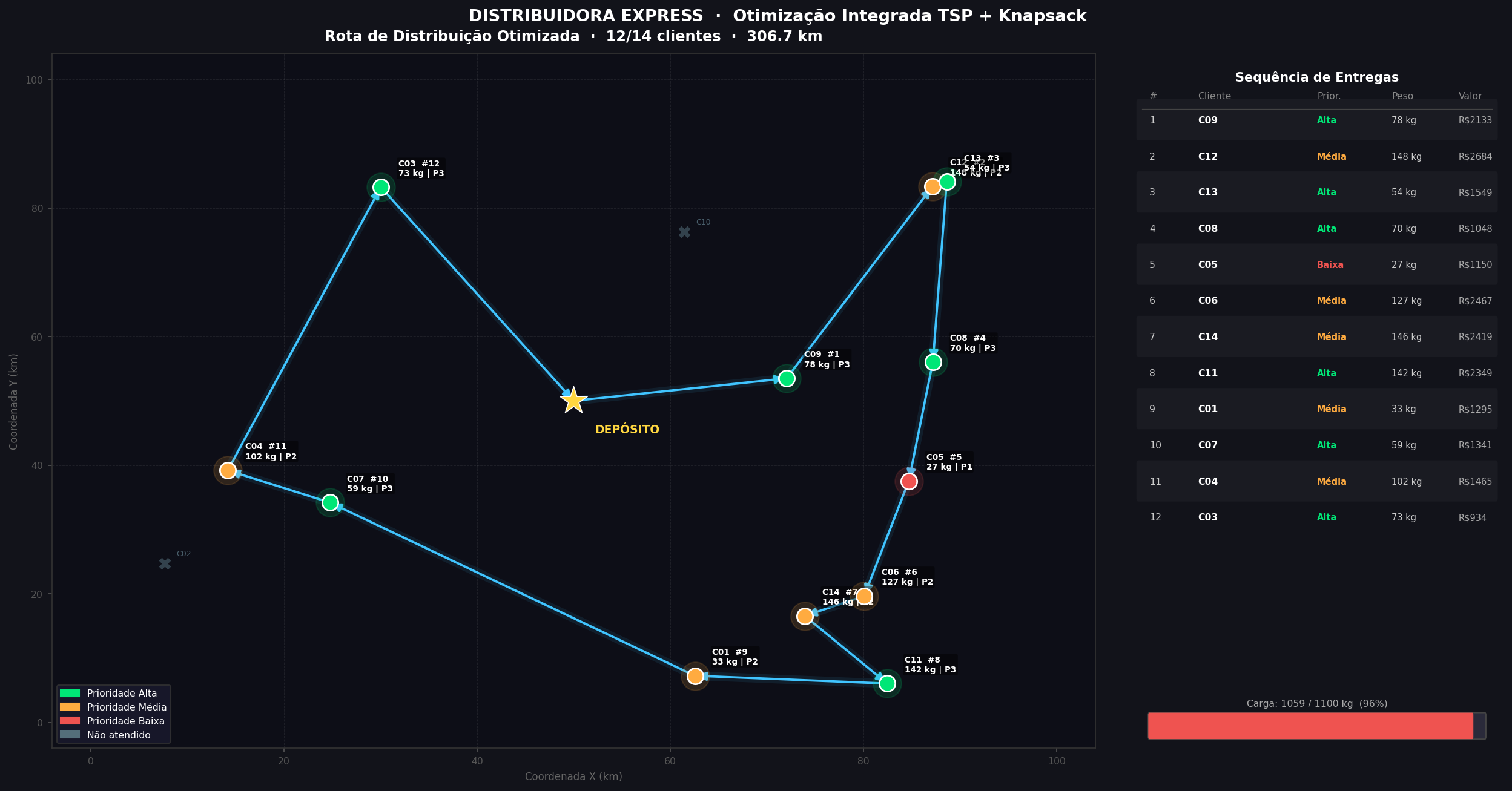}
  \caption{Logistics map: route, served customers (color = priority), excluded customers (gray X).}
  \label{fig:mapa}
\end{figure}

\begin{figure}[H]
  \centering
  \includegraphics[width=0.98\textwidth]{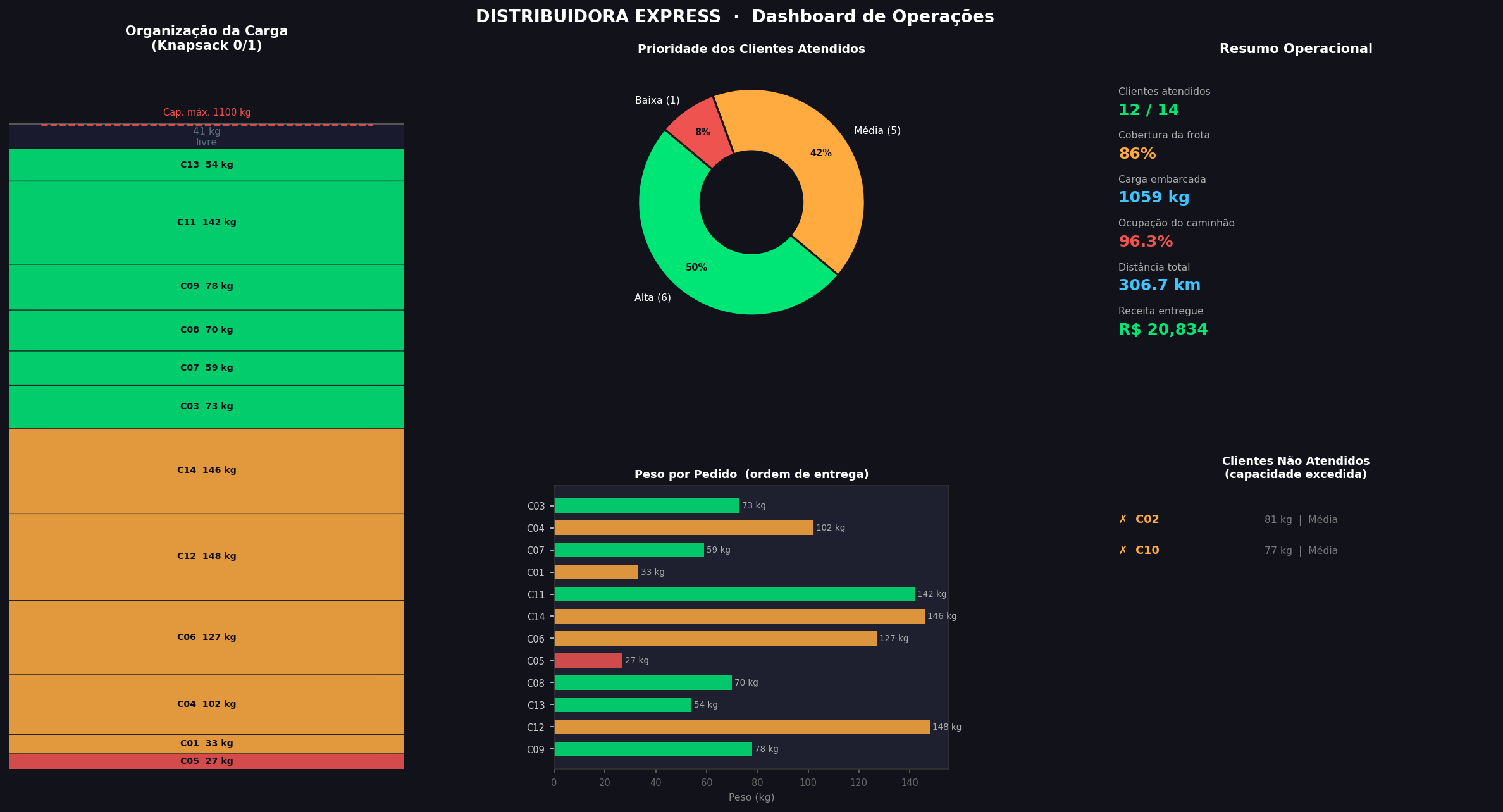}
  \caption{Dashboard: (A) cargo utilization; (B) priority distribution; (C) weight per order along the sequence; (D) KPIs; (E) unserved customers.}
  \label{fig:dashboard}
\end{figure}

Figures~\ref{fig:mapa} and~\ref{fig:dashboard} present the logistics map and the operational dashboard. Visual details are described briefly; the quantitative analysis is in the tables.

Panel~A of Figure~\ref{fig:dashboard} represents the load as a vertical column: the 40~kg slack (3.6\% of capacity) highlights the efficiency of the DP. Panel~B confirms the dominance of medium and high priorities among served customers. Panel~D summarizes the KPIs: coverage $9/14 \approx 64\%$, utilization $1060/1100 \approx 96\%$, revenue R\$~18,392.60.

\section{Analysis of Parameter $\alpha$}

Table~\ref{tab:alpha} analyzes the effect of $\alpha$ on the VNN-P score for two customers at 10~km, with priorities 3 and 1.

\begin{table}[H]
\centering
\caption{Effect of $\alpha$ on the score (distance = 10~km)}
\label{tab:alpha}
\begin{tabular}{@{}cccccl@{}}
\toprule
$\alpha$ & $s_{ij}$ ($p=3$) & $s_{ij}$ ($p=1$) & Factor & Regime \\
\midrule
0.00 & 10.00 & 10.00 & $1.00\times$ & Classical VNN \\
0.25 & 7.60  & 10.00 & $1.32\times$ & Mild bias \\
0.55 & 5.25  & 10.00 & $1.90\times$ & Balance (adopted) \\
1.00 & 3.33  & 10.00 & $3.00\times$ & Strong bias \\
2.00 & 1.11  & 10.00 & $9.00\times$ & Quasi-priority-first \\
\bottomrule
\end{tabular}
\end{table}

For $\alpha=0.55$, a high-priority customer is preferred if it is up to $3^{0.55} \approx 1.90$ times farther than a low-priority one. This defines an interpretable \textbf{geographic indifference zone}. A systematic analysis of the effect of $\alpha$ on total distance and visiting order (via experiments in the protocol of Section~\ref{sec:protocolo}) is left for future work.

\section{Conclusion, Limitations, and Future Work}

This work presented an applied routing variant with customer selection, combining priority-aware 0/1 Knapsack with priority-biased TSP. The contributions are: (i) quadratic priority weighting in the Knapsack objective function, with sensitivity analysis; (ii) a parameterizable priority score in the VNN-P, with geometric interpretation; (iii) a comparative experimental protocol that allows assessing the quality of the hierarchical decomposition; and (iv) an explicit discussion of the limitations.

\noindent\textbf{Limitations.} The hierarchical decomposition (Knapsack $\to$ TSP) does not guarantee global optimality: a low-priority customer close to the route may be excluded, increasing total distance. Moreover, the base scenario is single and small-scale; generalization to larger instances requires the experiments proposed in Section~\ref{sec:protocolo}. The parameter $\alpha$ was calibrated heuristically, without systematic search.

\noindent\textbf{Future work.} (i) Execute the full experimental protocol and report comparative results; (ii) extend to multiple vehicles (CVRP) and time windows (VRPTW); (iii) apply metaheuristics for PTP/OP \citep{Chagas2021} to the TSP component on larger instances; (iv) use public OP benchmarks \citep{Vansteenwegen2011,Gunawan2016} for comparison; (v) integrate real geospatial data via OpenStreetMap.

\bibliographystyle{plainnat}
\bibliography{references}

\end{document}